\documentclass[leqno,12pt]{article}
\usepackage{amssymb,amsfonts}
\usepackage{amsmath,latexsym}

\usepackage{amstext,epic,eepic,epsf,pslatex}
\usepackage{graphicx}

\newcommand{\bepr}{{\em Proof} } 
\newcommand{\enpr}{\hfill \rule{.5em}{.5em}}

\newcommand{\R}{{\mathbb R}}

\newcommand{\Z}{{\mathbb Z}}

\newcommand{\Tr}{\hbox{Tr\,}}

\def\Xint#1{\mathchoice
{\XXint\displaystyle\textstyle{#1}}%
{\XXint\textstyle\scriptstyle{#1}}%
{\XXint\scriptstyle\scriptscriptstyle{#1}}%
{\XXint\scriptscriptstyle\scriptscriptstyle{#1}}%
\!\int}
\def\XXint#1#2#3{{\setbox0=\hbox{$#1{#2#3}{\int}$ }
\vcenter{\hbox{$#2#3$ }}\kern-.6\wd0}}

\newtheorem{defin}{Definition}[section] 
\newtheorem{prop}{Proposition}[section] 
\newtheorem{thm}{Theorem}[section] 
\newtheorem{lemma}{Lemma}[section]

\newtheorem{open}{Open Question}[section] 

\newtheorem{cor}{Corollary}[section]

\begin{document}

\title{Divergence-free positive symmetric tensors and fluid dynamics}

\author{Denis Serre \\ \'Ecole Normale Sup\'erieure de Lyon\thanks{U.M.P.A., UMR CNRS--ENSL \# 5669. 46 all\'ee d'Italie, 69364 Lyon cedex 07. France. {\tt denis.serre@ens-lyon.fr}}}

\date{\today}

\maketitle

\begin{abstract}
We consider $d\times d$ tensors $A(x)$ that are symmetric, positive semi-definite, and whose row-divergence vanishes identically. We establish sharp inequalities for the integral of $(\det A)^{\frac1{d-1}}$. We apply them to models of compressible inviscid fluids: Euler equations, Euler--Fourier, relativistic Euler, Boltzman, BGK, etc... We deduce an {\em a priori} estimate for a new quantity, namely the space-time integral of $\rho^{\frac1n}p$, where $\rho$ is the mass density, $p$ the pressure and $n$ the space dimension. For kinetic models, the corresponding quantity generalizes Bony's functional.
\end{abstract}


\paragraph{Notations.}
The integer $d\ge2$ is the number of independent variables, which are often space-time coordinates. It serves also for the size of square matrices. If $1\le j\le d$ and $x\in\R^d$ are given, we set $\widehat{x_j}=(\ldots,x_{j-1},x_{j+1},\ldots)\in\R^{d-1}$~; the projection $x\mapsto\widehat{x_j}$ ignores the $j$-th coordinate. The transpose of a matrix $M$ is $M^T$.  If $A\in {\bf M}_d(\R)$, its cofactor matrix $\widehat{A}$ satisfies
$$\widehat{A}^TA=A\,\widehat{A}^T=(\det A)I_d,\qquad\det\widehat{A}=(\det A)^{d-1}.$$
Because we shall deal only with symmetric matrices, we have simply $\widehat{A}\,A=A\,\widehat{A}=(\det A)I_d$.
The space of $d\times d$ symmetric matrices with real entries is ${\bf Sym}_d$. The cones of positive definite, respectively positive semi-definite, matrices are  ${\bf SPD}_d$ and ${\bf Sym}_d^+$. If $u\in\R^d$, $u\otimes u\in{\bf Sym}_d^+$ denotes the rank-one matrix of entries $u_iu_j$.

The unit sphere of $\R^d$ is $S^{d-1}$. The Euclidian volume of  an open subset $\Omega$ of $\R^d$ is denoted $|\Omega|$. If the boundary $\partial\Omega$ is rectifiable, we denote the same way $|\partial\Omega|$ its area, and $ds(x)$ the area element. For instance, the ball $B_r$ of radius $r$ and its boundary, the sphere $S_r$, satisfy $|B_r|=\frac rd\,|S_r|$. If $\Omega$ has a Lipschitz boundary, its outer unit normal $\vec n$ is defined almost everywhere.

If $f:\Omega\rightarrow\R$ is integrable, its average over $\Omega$ is the number
$$\Xint-_\Omega f(x)\,dx:=\frac1{|\Omega|}\,\int_\Omega f(x)\,dx.$$
Given a lattice $\Gamma$ of $\R^d$, and $f:\R^d\rightarrow\R$ a $\Gamma$-periodic, locally integrable function, we denote
$$\int_{\R^d/\Gamma}f(x)\,dx$$
the value of the integral of $f$ over any fundamental domain. We define as above the average value
$$\Xint-_{\R^d/\Gamma}f(x)\,dx.$$

For our purpose, a {\em tensor} is a matrix-valued function $x\mapsto T(x)\in{\bf M}_{p\times q}(\R)$. If $q=d$ and if the derivatives make sense (say as distributions), we form
$${\rm Div}T=\left(\sum_{j=1}^d\partial_jt_{ij}\right)_{1\le i\le p},$$
which is vector-valued.
We emphasize the uppercase letter D in this context. We reserve the lower case operator {\bf div} for vector fields.

If $1\le p\le\infty$, its conjugate exponent is $p'$.

\paragraph{Acknowledgements.} This research benefited from discussions that I had with various persons, and references that I got from others. I thank warmly Gr\'egoire Allaire, Yann Brenier, Vincent Calvez, Guido de Philippis, Reinhard Illner, Gr\'egoire Loeper, Petru Mironescu, Jean-Christophe Mourrat, Laure Saint-Raymond, Bruno S\'evennec and C\'edric Villani.

\section{Motivations}

We first define the mathematical object under consideration.
\begin{defin} Let $\Omega$ be an open subset of $\R^d$. A {\em divergence-free positive symmetric tensor} (in short, a {\em DPT}) is a locally integrable tensor $x\mapsto A(x)$ over $\Omega$ with the properties that $A(x)\in{\bf Sym}_d^+$ almost everywhere, and ${\rm Div}\,A=0$.
\end{defin}
The following fact is obvious.
\begin{lemma}[Congruence.]
If $A$ is a DPT and $P\in{\bf GL}_d(\R)$ is given, then the tensor $$B(y):=PA(P^{-1}y)P^T$$ is also a DPT.
\end{lemma}

\subsection{Where do the divergence-free positive symmetric tensors occur ?}

Most of our examples, though not all of them, come from fluid dynamics, where a DPT contains a stress tensor.
\begin{description}
\item[Compressible gas.] In space dimension $n\ge1$, a gas is described by a mass density $\rho\ge0$, a velocity $u$ and a pressure $p\ge0$. These fields obey the Euler equations (conservation of mass and momentum)
$$\partial_t\rho+{\rm div}_y(\rho u)=0,\qquad\partial_t(\rho u)+{\rm Div}_y(\rho u\otimes u)+\nabla_y p=0.$$
 Here $x=(t,y)$ and $d=1+n$. The tensor
$$A(t,y)=\begin{pmatrix} \rho & \rho u^T \\ \rho u & \rho u\otimes u +pI_n \end{pmatrix}$$
is a DPT.
\item[Rarefied gas.] It is described by a density function $f(t,y,v)\ge0$ where $v\in\R^n$ is the particle velocity. The evolution is governed by a kinetic equation
$$(\partial_t+v\cdot\nabla_y)f=Q[f(t,y,\cdot)].$$
The left-hand side is the transport operator, while the right-hand side, a non-local operator acting on the velocity variable, accounts for the interaction between particles. This class  contains the Boltzman equation, as well as the discrete kinetic models or the BGK model. When the collisions are elastic, the mass, momentum and energy are conserved. This is reflected by the properties
$$\int_{\R^n}Q[g](v)\,dv=0,\qquad\int_{\R^n}Q[g](v)v\,dv=0,\qquad\int_{\R^n}Q[g](v)|v|^2\,dv=0$$
for every reasonable function $g(v)$.
Integrating the kinetic equation against $dv$, $v\,dv$ and $\frac12|v|^2dv$, we obtain again, at least formally, the conservation laws
$$\partial_t\rho+{\rm div}_ym=0,\qquad\partial_tm+{\rm Div}_yT=0,\qquad\partial_tE+{\rm div}_yQ=0,$$
where 
$$\rho(t,y):=\int_{\R^n}f(t,x,v)\,dv,\qquad m(t,y):=\int_{\R^n}f(t,x,v)\,v\,dv,\qquad E:=\int_{\R^n}f(t,y,v)\,\frac12|v|^2dv$$
are the mass density, linear momentum and energy, while 
$$T:=\int_{\R^n}f(t,y,v)\,v\otimes v\,dv,\qquad Q:=\int_{\R^n}f(t,y,v)\,\frac12|v|^2v\,dv$$
are fluxes. The tensor
$$A(t,y)=\begin{pmatrix} \rho & m^T \\ m & T \end{pmatrix}$$
is again a DPT.
\item[Steady\,/\,self-similar flows.] Let us go back to gas dynamics. If the flow is steady, then on the one hand ${\rm div}(\rho u)=0$, and on the other hand ${\rm Div}(\rho u\otimes u)+\nabla p=0$. Therefore the tensor $A=\rho u\otimes u+pI_n$ is a DPT in the physical domain $\Omega\subset\R^n$.

If instead the flow is self-similar, in the sense that $\rho,u$ and $p$ depend only upon $\xi=\frac yt$ (this is reminiscent to the multi-D Riemann Problem), then it obeys to the reduced system
\begin{equation}\label{RP}
{\rm div}_\xi(\rho v)+n\rho=0,\qquad {\rm Div}_\xi(\rho v\otimes v)+\nabla_\xi p+(n+1)\rho v=0,
\end{equation}
where $v:=u(\xi)-\xi$ is the {\em pseudo-velocity}. The tensor $A:=\rho v\otimes v+pI_n$ is not a DPT, because of the source term $(n+1)\rho v$. However it is positive semi-definite, and we shall be able to handle such a situation.
\item[Relastivistic gas dynamics.]
In the Minkowski space, the Euler equations write ${\rm Div}\,T=0$ where $T$ is the stress-energy tensor. This is another instance of a DPT.
\end{description}

\paragraph{Periodic homogenization of elliptic operators.}
This is a completely different context, for which we refer to \cite{Gal,Tar}. A $\Gamma$-periodic symmetric tensor $A(x)$ is given, which satisfies the bounds
$$\alpha|\xi|^2\le\xi^TA(x)\xi\le\beta|\xi|^2,\qquad\forall\,\xi\in\R^d,$$
where $0<\alpha\le\beta<+\infty$ are constants.
The differential operator $Lu={\rm div}(A\nabla u)$ is uniformly elliptic. Given a vector $\xi$, the problem
$${\rm div}(A(\xi+\nabla u))=0$$
admits a unique $\Gamma$-periodic solution $u_\xi\in H^1_{loc}$, up to an additive constant. 
A PDE such as (\ref{elpsi}) below governs the temperature or the electric potential at equilibrium in a periodic non-homogeneous medium. The macroscopic behaviour of the medium is well described by the so-called {\em effective tensor} $A_{\rm eff}$, whose definition is
$$A_{\rm eff}\xi=\Xint-_{\R^d/\Gamma}A(x)(\xi+\nabla u_\xi)\,dx.$$
An equivalent formulation is
\begin{equation}\label{Aeffxi}
\xi^TA_{\rm eff}\xi=\Xint-_{\R^d/\Gamma}(\xi+\nabla u_\xi)^TA(x)(\xi+\nabla u_\xi)\,dx=\inf_{w\in H^1_{per}}\Xint-_{\R^d/\Gamma}(\xi+\nabla w)^TA(x)(\xi+\nabla w)\,dx.
\end{equation}
In particular, $A_{\rm eff}\in{\bf SPD}_d$.  The effective tensor is known to obey  the sharp bounds
\begin{equation}\label{boundsA}
A_-\le A_{\rm eff}\le A_+
\end{equation}
where $A_\pm$ are the harmonic and arithmetic means of $A(x)$~:
$$A_+=\Xint-_{\R^d/\Gamma}A(x)\,dx,\qquad A_-=\left(\Xint-_{\R^d/\Gamma}A(x)^{-1}\,dx\right)^{-1}.$$
\begin{prop}
The effective tensor $A_{\rm eff}$ equals the upper bound $A_+$ if, and only if, $A$ is a DPT.
\end{prop}

Although this is a classical and simple fact, we recall the proof.
Taking $w\equiv0$ in (\ref{Aeffxi}), we obtain the upper bound $\xi^TA_{\rm eff}\xi\le\xi^TA_+\xi$. If $A_{\rm eff}=A_+$, this implies that the infimum is attained precisely at constants~; in other words $\nabla u_\xi\equiv0$. But then ${\rm div}(A(\xi+\nabla u_\xi))=0$ writes ${\rm div}(A\xi)=0$. This being true for every $\xi$, we have ${\rm Div}\,A=0$. The converse is immediate: if ${\rm Div}\,A=0$, then $u_\xi$ is just a constant, and therefore $\xi^TA_{\rm eff}\xi=\xi^TA_+\xi$.

\bigskip

The role of the effective tensor is the following. Given $f\in H^{-1}(\Omega)$ and a small scale $\epsilon>0$, the solution $u^\epsilon$ of the Dirichlet boundary-value problem
\begin{equation}\label{elpsi}
{\rm div}\left(A(\frac x\epsilon)\nabla u^\epsilon\right)=f(x),\qquad u^\epsilon|_{\partial\Omega}=0
\end{equation}
remains bounded in $H^1(\Omega)$ and converges weakly as $\epsilon\rightarrow0$ towards the solution $\bar u$ of the same problem with the effective matrix:
$${\rm div}\left(A_{\rm eff}\nabla \bar u\right)=f(x),\qquad \bar u|_{\partial\Omega}=0.$$
When $f\in L^2(\Omega)$ instead, the sequence $u^\epsilon$ remains bounded in $H^2(\Omega)$ only if $A_{\rm eff}$ coincides with $A_+$, see \cite{CV}. This is due to the fact that the first corrector in the expansion of $u^\epsilon$ in terms of $\epsilon$ becomes trivial.

\subsubsection{When divergence-free symmetric tensors are not positive}

It is fair to list a few important examples in which our approach does not apply because of the lack of positiveness.

\paragraph{Compressible Navier-Stokes equations.}
The system that governs a viscous compressible fluid differs slightly from the Euler equation. The conservation of mass remains the same, but the conservation of momentum becomes
$$\partial_t(\rho u)+{\rm Div}_y(\rho u\otimes u-\lambda(\nabla u^T+\nabla u))+\nabla_y (p-\mu\,{\rm div}\,u)=0.$$
The divergence-free tensor
$$\begin{pmatrix} \rho & \rho u^T \\ \rho u & \rho u\otimes u-\lambda(\nabla u^T+\nabla u)+(p-\mu\,{\rm div}\,u)\,I_n \end{pmatrix}$$
is not positive in general.

\paragraph{Mean-field equations.}
One form of kinetic models is
\begin{equation}\label{MF}
(\partial_t+v\cdot\nabla_y)f+F(t,y)\cdot\nabla_vf=0,
\end{equation}
where the force $F$ is coupled to the density $\rho=\int f\,dv$ through $F=-\nabla_yE$,
$$E:=\phi*_y\rho=\int_{\R^n}\phi(y-z)\rho(z)\,dz.$$
The potential $\phi$ is a characteristic of the model. For instance a Coulomb force or the gravity yield the coupling
$$\Delta E=\beta\rho$$
where $\beta$ is a constant that can be positive (attractive force) or negative (repulsive force). With this choice, (\ref{MF}) implies formally the hydrodynamic system
$$\partial_t\rho+{\rm div}_ym=0,\qquad\partial_tm+{\rm Div}_yT=0,$$
where as usual $\rho$ and $m$ are the moments of $f$ of order $0$ and $1$, and 
$$T=\int_{\R^n}f(t,y,v)v\otimes v\,dv+\frac1\beta\,(F\otimes F-\frac12\,|F|^2I_n).$$
Because $T$ does not have a definite sign, the tensor $\begin{pmatrix} \rho & m^T \\ m & T \end{pmatrix}$ is not positive in general.

\paragraph{Added in proofs.} Here is a short list of divergence-free symmetric tensors in other models from physics or mechanics. The energy-momentum tensor of the electromagnetic field in vaccum, when normalizing the light speed to $c=1$~; its symmetry is related to the Lorentz invariance of the Lagrangian $\omega\mapsto L(\omega)$ where $\omega=(E\cdot dx)\times dt+(B\times dx)\cdot dx$ denotes the electromagnetic field. The Lagrangian needs not be quadratic. The mass-momentum tensor in a Schr\"odinger equation. The energy-momentum tensor in hyper-elasticity, written in Eulerian coordinates~; the symmetry is related to the conservation of angular momentum (frame indifference). Only the last one may be positive semi-definite~; this arises when the stored energy $\varepsilon(F^TF)$ ($F$ the deformation tensor) is a monotonous non-increasing function of $C:=F^TF$. This usually requires that the medium be compressed, $C\le I_3$.

\subsection{$\Lambda$-concave functions}

Let $K$ be a convex subset of some space $\R^N$ and $F:K\rightarrow\R$ be a continuous function. We consider measurable functions $u:\Omega\rightarrow K$ (say, bounded ones). Let us recall that $F$ is concave if, and only if the inequality
\begin{equation}\label{FmoyF}
\Xint-_\Omega F(u)\,dx\le F\left(\Xint-_\Omega u\,dx\right)
\end{equation}
for every such $u$. This is just a reformulation of Jensen's inequality. In particular, the equality holds true for every $u$ if, and only if $F$ is affine.

A general question, first addressed by F. Murat and L. Tartar \cite{Mur,Tar_HW} is whether a differential constraint imposed to $u$ allows some non-concave functions $F$ to satisfy (\ref{FmoyF}). For instance, the following is known \cite{Ball}. If $\Omega=\R^d/\Gamma$, and $u=\nabla\phi$ (hence $F$ applies to $d\times m$ matrices, and ${\rm curl}\,u=0$) is $\Gamma$-periodic, then the equality holds true in (\ref{FmoyF}) whenever $F$ is a linear combination of minors. And the inequality is valid for every {\em polyconcave} function, that is a concave function of all the minors.

The same question is addressed here, when $\R^N={\bf Sym}_d$, the cone $K$ is ${\bf Sym}_d^+$ and the differential constraint is ${\rm Div}\,A=0$.  Every concave function satisfies it, in a trivial manner because the inequality does not involve the differential constraint. A fundamental example of that situation is the function
$$A\mapsto(\det A)^{\frac1d},$$
which is concave over ${\bf Sym}_d^+$ (see  \cite{Ser_mat} Section 6.6).

\paragraph{A necessary condition.}
Let us recall a construction due to Tartar \cite{Tar_HW}. Let $B,C\in{\bf Sym}_d^+$ be given, such that $C-B$ is singular (that is $\det(C-B)=0$). Then there exists a non-zero vector $\xi$ such that $(C-B)\xi=0$. This ensures that  for every function $g:\R\rightarrow\{0,1\}$, the tensor
$$A(x):=g(x\cdot\xi)B+(1-g(x\cdot\xi))C$$
is a DPT. If $F$ satisfies (\ref{FmoyF}) then in particular we have
$$\Xint-_{\R^d/\Z^d}F(g(x\cdot\xi)B+(1-g(x\cdot\xi))C)\,dx\le F\left(\Xint-_{\R^d/\Z^d}(g(x\cdot\xi)B+(1-g(x\cdot\xi))C)\,dx\right).$$
With $\theta$ the mean value of $g$, this is
$$\theta F(B)+(1-\theta)F(C)\le F(\theta B+(1-\theta)C).$$
The restriction of $F$ to the segment $[B,C]$ must therefore be concave. We say that $F$ is {\em $\Lambda$-concave}, where $\Lambda$ is the cone of singular symmetric matrices.

\bigskip

Let us go back to the trivial example of $A\mapsto(\det A)^{\frac1d}$. Is it possible to improve the exponent $\frac1d$ while keeping the $\Lambda$-concavity~? The answer is positive:
\begin{prop}\label{p:Fal}
For an exponent $\alpha>0$, the map
\begin{eqnarray*}
A & \longmapsto & (\det A)^\alpha \\
{\bf Sym}_d^+ & \rightarrow & \R^+
\end{eqnarray*}
is $\Lambda$-concave if, and only if $\alpha\le\frac1{d-1}$\,.
\end{prop}

\bepr

Let $A,A+B\in{\bf Sym}_d^+$ be such that $\det B=0$ and denote $f(t)=(\det(A+tB))^{\frac1{d-1}}$.
To prove that $f$ is concave over $[0,1]$, it is enough to prove that $f(t)\le f(0)+tf'(0)$. Using a congruence, we may assume that $A=I_d$. Another congruence, by an orthogonal matrix $P$, allows us to assume that in addition, $B$ is diagonal: $B={\rm diag}(b_1,\ldots,b_{d-1},0)$. Then, using the geometric-arithmetic mean inequality,
$$f(t)=\prod_{j=1}^{d-1}(1+tb_j)^{\frac1{d-1}}\le \frac1{d-1}\sum_{j=1}^{d-1}(1+tb_j)=f(0)+tf'(0).$$

If $\alpha<\frac1{d-1}$\,, then the function $F_\alpha$ under consideration is a composition $\phi_\alpha\circ F_{\frac1{d-1}}$ where $\phi_\alpha(s)=s^{\frac\alpha{d-1}}$. Since $\phi_\alpha$ is concave increasing and $F_{\frac1{d-1}}$ is concave, $F_\alpha$ is concave. 

Conversely, if $F_\alpha$ is $\Lambda$-concave and $B={\rm diag}(b_1,\ldots,b_{d-1},0)$ is singular, diagonal with all $b_j>0$, then 
$$t\mapsto\prod_{j=1}^{d-1}(1+tb_j)^\alpha$$
must be concave. In particular it must be sub-linear, which implies $\alpha\le\frac1{d-1}$\,.

\enpr

\bigskip

Once we know that $F_\alpha$ passes the test of $\Lambda$-concavity, it becomes natural to ask whether it satisfies a functional inequality, such as (\ref{FmoyF}) when $\Omega=\R^d/\Gamma$, or something similar when $\Omega$ is a bounded domain.

Clues  are provided by two particular cases:
\begin{description}
\item[Diagonal case.] A diagonal DPT is a map $x\mapsto{\rm diag}(g_1(\widehat{x_1}),\ldots,g_d(\widehat{x_d}))$, where the $j$-th function (non-negative) does not depend upon $x_j$. Such a tensor is periodic whenever the $g_j$'s are so, and the lattice is parallel to the axes. This situation enjoys an inequality due to Gagliardo \cite{Gag}:
\begin{equation}\label{Gagl}
\Xint-_{\R^d/\Gamma}\left(\prod_{j=1}^dg_j(\widehat{x_j})\right)^{\frac1{d-1}}dx\le\prod_{j=1}^d\left(\Xint-_{\R^{d-1}/\Gamma_j}g_j(\widehat{x_j})\,d\widehat{x_j}\right)^{\frac1{d-1}},
\end{equation}
where the lattice $\Gamma_j$ is the projection of $\Gamma$ on the hyperplane $x_j=0$. The right-hand side can be viewed as the average of a power of $\det A$, while the left-hand side is the power of the determinant of the average matrix.
\item[Cofactors of Hessian.] Let $\phi\in W^{2,d-1}(\Omega)$ be a convex function over a convex domain $\Omega$. Let us form its Hessian matrix $\nabla^2\phi$, and then the positive symmetric tensor $A=\widehat{\nabla^2\phi}$.
\begin{lemma}\label{l:widehat} The tensor defined above is a DPT.
\end{lemma}
The proof consists in remarking that the differential form $\omega_j:=\sum_ia_{ij}dx_j$ is nothing but the exterior product $\cdots\wedge d\phi_{j-1}\wedge d\phi_{j+1}\wedge\cdots$, where only the factor $d\phi_j$ has been omitted. This $(d-1)$-form is obviously closed, and this translates into the identity $\sum_i\partial_ia_{ij}=0$.

It turns out that $(\det A)^{\frac1{d-1}}=\det\nabla^2\phi$ is itself an exterior derivative, for instance that of $\phi_j\omega_j$. Therefore
$$\int_\Omega(\det A)^{\frac1{d-1}}dx$$
is actually a boundary integral. 

In the periodic case, we assume that only $\nabla^2\phi$ is $\Gamma$-periodic, and we write $\phi(x)=\frac12 x^TSx+\hbox{linear}+\psi(x)$ where $\psi$ is $\Gamma$-periodic. Then we have
$$\Xint-_{\R^d/\Gamma}(\det A)^{\frac1{d-1}}dx = \Xint-_{\R^d/\Gamma}\det(S+\nabla^2\psi)\,dx =\det S, $$
because the determinant of $S+\nabla\phi^2$ is the sum of $\det S$ and a linear combination of minors of $\nabla^2\phi$, each one being a divergence, thus integrating to zero. On the other hand we have
$$\Xint-_{\R^d/\Gamma}A(x)\,dx=\Xint-_{\R^d/\Gamma}\widehat{S+\nabla^2\psi}=\widehat{S}$$
for the same reason. We infer a remarkable identity:
\begin{prop}\label{p:idHess}
The formula $A=\widehat{S+\nabla^2\psi}$, where $\psi$ is $\Gamma$-periodic and $x\mapsto\frac12\,x^TSx+\psi(x)$ is convex, provides a DPT, which satisfies
$$\Xint-_{\R^d/\Gamma}(\det A)^{\frac1{d-1}}dx =\left(\det\Xint-_{\R^d/\Gamma} A(x)\,dx\right)^{\frac1{d-1}}.$$
\end{prop}
\end{description}

Both particular cases above are given in a periodic context but have counterparts in bounded convex domains. We shall explain below how they embed into results that are valid for every DPT.
The version in a bounded convex domain will involve the trace $A\vec n$, an object that makes sense just because of the divergence-free assumption.

\bigskip

The next two sections contain our results. Up to our knowledge, they have not been uncovered so far, perhaps because the DPT structure has been overlooked, or has been examined only at the linear level. Our results are two-fold. On the one hand we make general statements about DPTs, which are proved in Sections \ref{s:convex}  and \ref{s:per}.  The moral of these results is that the row-wise divergence operator displays a small amount of ellipticity~; when a control of ${\rm Div}\,A$ is coupled with the assumption of symmetry and semi-definiteness, then $\det A$ enjoys a slightly better integrability than $A$ itself. On the other hand, we give several applications to gas dynamics. They concern either the Euler system of a compressible fluid, or the kinetic models, for instance that of Boltzmann. Details are given in Section \ref{s:fluid}.

\section{General statements}

We present two abstract results about DPTs, which cover the periodic case and that of a convex bounded domain. The central object here is the application $F_{\frac1{d-1}}$~:
\begin{eqnarray*}
A & \longmapsto & (\det A)^{\frac1{d-1}}, \\
{\bf Sym}_d^+ & \longrightarrow & \R^+
\end{eqnarray*}

\subsection{Periodic case}

\begin{thm}\label{th:per}
Let the DPT $x\longmapsto A(x)$ be $\Gamma$-periodic, with $A\in L^1(\R^d/\Gamma)$. Then $(\det A)^{\frac1{d-1}}\in L^1(\R^d/\Gamma)$ and there holds
\begin{equation}\label{Jdetper}
\Xint-_{\R^d/\Gamma}(\det A(x))^{\frac1{d-1}}dx\le\left(\det\Xint-_{\R^d/\Gamma}A(x)\,dx\right)^{\frac1{d-1}}.
\end{equation}
\end{thm}

\bigskip

An easy consequence is the following, which displays a little gain of integrability.
\begin{cor}\label{c:CS}
Let $\Omega$ be an open set of $\R^d$. Let $\bar A\in{\bf SPD}_d$ be given, and  $A$ be a DPT over $\Omega$, such that $A-\bar A$ is compactly supported. Then 
$$\Xint-_\Omega A(x)\,dx=\bar A\quad\hbox{and}\quad\Xint-_{\Omega}(\det A(x))^{\frac1{d-1}}dx\le\left(\det\Xint-_{\Omega}A(x)\,dx\right)^{\frac1{d-1}}.$$
\end{cor}
The inequality (\ref{Jdetper}) of Theorem \ref{th:per} is actually sharp:
\begin{prop}\label{p:caseeq}
In the situation of Theorem \ref{th:per}, suppose that $x\mapsto\det A$ is a smooth function, bounded by below and by above. Then the equality case in (\ref{Jdetper}) is achieved if, and only if $A=\widehat{\nabla^2\theta}$, where $\theta$ is a convex function whose Hessian is periodic.
\end{prop}
We expect that the assumptions that $\det A$ is smooth and bounded below by a positive constant can be removed, though we do not dwell into more details here. 

\bigskip

Another interesting consequence is the following (recall that $d\ge2$).
\begin{cor}\label{c:surprise}
Let $\theta\in W^{2,d-1}_{\rm loc}(\R^d)$ be a convex function, whose Hessian is $\Gamma$-periodic. Then $\det{\rm D}^2\theta$ is integrable over $\R^d/\Gamma$.
\end{cor}
{\em Proof}: Just apply Theorem \ref{th:per} to $A=\widehat{{\rm D}^2\theta}$.

\bigskip

Within the context of periodic homogenization, (\ref{Jdetper}) applies to the case where $A_{\rm eff}=A_+$. One might wonder whether it is a particular case of a more general inequality, once $A_{\rm eff}$ differs from $A_+$. We leave this question open, but it is easy to rule out the tempting inequality
\begin{equation}\label{tempt}
\Xint-_{\R^d/\Gamma}(\det A(x))^{\frac1{d-1}}dx\quad\stackrel?\le\quad\left(\det A_{\rm eff}\right)^{\frac1{d-1}}.
\end{equation}
As a matter of fact, the upper bound in (\ref{boundsA}) and the monotonicity of the determinant tell us that 
$$\det A_{\rm eff}\le\det\Xint-_{\R^d/\Gamma}A(x)\,dx.$$
If the inequality (\ref{tempt}) was true, then the average of $F(A):=(\det A)^{\frac1{d-1}}$ would be bounded above by $F$ of the average, for every $x\mapsto A(x)$ taking values in ${\bf SPD}_d$. This would imply the concavity of $F=F_{\frac1{d-1}}$ over ${\bf SPD}_d$, which we know is false (Proposition \ref{p:Fal}).

\bigskip

\subsubsection{Extension to general symmetric positive tensors}

When $A$ is not divergence-free, we still have the following surprising result.
\begin{thm}\label{th:nondiv}
Let $x\longmapsto A(x)$ be $\Gamma$-periodic, taking values in ${\bf Sym}_d^+$. Assume that $A\in L^1(\R^d/\Gamma)$ and ${\rm Div}\, A$ is a bounded measure over $\R^d/\Gamma$.
Then 
$$(\det A)^{\frac1{d-1}}\in L^1(\R^d/\Gamma).$$
\end{thm}

Theorem \ref{th:nondiv} can be compared with Sobolev embeddings and elliptic regularity. If the assumption that ${\rm Div}\,A$ is a bounded measure is replaced by the fact that every derivative $\partial_ia_{jk}$ is a bounded measure, then $A\in L^{\frac d{d-1}}$ and the conclusion follows immediately. Even if we only assume that $P(D)A$ is integrable for some elliptic first-order differential operator $P(D)$, we know that $A\in L^{\frac d{d-1}-\epsilon}$ for every $\epsilon>0$. The Theorem says that the operator {\bf Div} displays a (very weak) form of ellipticity, when combined to the symmetry and positivity of the tensor.

This comparison leads us to the following question, which we leave open.
\begin{quote}
\begin{open}
In Theorem \ref{th:nondiv}, assume instead that $A$ and ${\rm Div}\,A$ belong to $L^p(\R^d/\Gamma)$ with $1<p<d$. Is it true that $(\det A)^{\frac1d}$ belongs to $L^{p^*}(\R^d/\Gamma)$ with $\frac1{p^*}=\frac1p-\frac1d$\,?
\end{open}
\end{quote}

\subsection{Bounded domain}

We assume now that the domain $\Omega$ is convex. We recall that if a divergence-free vector field $\vec q$ belongs to $L^p(\Omega)$, then it admits a normal trace $\gamma_\nu\vec q$ which belongs to the Sobolev space $W^{-\frac1p,p}(\partial\Omega)$. It is defined by duality, by the formula
$$\langle \gamma_\nu\vec q,\gamma_0w\rangle=\int_\Omega\vec q\cdot\nabla w\,dx,\qquad\forall w\in W^{1,p'}(\Omega),$$
where $\gamma_0$ is the standard trace operator from $W^{1,p'}(\Omega)$ into $W^{\frac1p,p'}(\partial\Omega)$.

When $\vec q$ is a smooth field, $\gamma_\nu\vec q$ coincides with the pointwise normal trace $\vec q|_{\partial\Omega}\cdot\vec n$. We say that $\vec q$ has an {\em  integrable normal trace} if the distribution $\gamma_\nu\vec q$ coincides with an integrable function~; then we write $\vec q\cdot\vec n$ instead. For instance, and this is the case below, the row-wise trace $\gamma_\nu A$ of a DPT of class $L^d(\Omega)$ makes sense in $W^{-\frac1d,d}(\partial\Omega)$, and we denote this trace $A\vec n$ when it is integrable.
\begin{thm}\label{th:conv}
Let $\Omega$ be a bounded convex open subset in $\R^d$. Let $A$ be a DPT over $\Omega$ that belongs to $L^{\frac d{d-1}}_{loc}(\R^d)$ and has an integrable normal trace. Then there holds
\begin{equation}\label{detconv}
\int_{\Omega}(\det A(x))^{\frac1{d-1}}dx\le\frac1{d|S^{d-1}|^{\frac1{d-1}}}\,\|A\vec n\|_{L^1(\partial\Omega)}^{\frac d{d-1}}.
\end{equation}

If $A$ is only symmetric non-negative, but ${\rm Div}\,A$ is a bounded measure (therefore $A$ is not a DPT), then we have 
\begin{equation}\label{convmeas}
\int_{\Omega}(\det A(x))^{\frac1{d-1}}dx\le\frac1{d|S^{d-1}|^{\frac1{d-1}}}\,\left(\|A\vec n\|_{L^1(\partial\Omega)}+\|{\rm Div}\,A\|_{{\cal M}(\Omega)}\right)^{\frac d{d-1}},
\end{equation}
where the second norm is the total mass of the measure $|{\rm Div}\,A|$.
\end{thm}

\bigskip

The inequalities (\ref{Jdetper}) and (\ref{detconv}) can be viewed as {\em non-commutative} analogues of the Gagliardo inequality (\ref{Gagl}).

\bigskip

Remark that a somehow more elegant form of (\ref{detconv}) happens when $\Omega$ is a ball:
\begin{equation}\label{elegant}
\Xint-_{B_r}(\det A(x))^{\frac1{d-1}}dx\le\left(\Xint-_{S_r}|A\vec n|\,ds(x)\right)^{\frac d{d-1}}.
\end{equation}

\bigskip

Once again, the inequality (\ref{detconv}) is sharp, and we have
\begin{prop}\label{p:eqcase}
In the situation of Theorem \ref{th:conv}, suppose that $x\mapsto\det A$ is a smooth function, bounded by below and by above. Then the equality case in (\ref{detconv}) is achieved if, and only if $A=\widehat{\nabla^2\theta}$, where $\theta$ is a convex function such that $\nabla\theta(\Omega)$ is a ball centered at the origin.
\end{prop}

\bigskip

On a qualitative side, we have the following result.
\begin{prop}\label{p:suppK}
Let $\Omega$ be a bounded open subset of $\R^d$ with a Lipschitz boundary. Let $A$ be a DPT over $\Omega$. If $\vec n^TA\vec n\equiv0$ over $\partial\Omega$, then $A$ vanishes identically over $\Omega$.
\end{prop}

\subsubsection{Gain of integrability}

The following result is more in the spirit of Theorem \ref{th:per}.
\begin{thm}\label{th:gain}
Let $\Omega$ be an open domain of $\R^d$, and $A$ be a symmetric, positive semi-definite tensor of class $L^1_{\rm loc}(\Omega)$ and such that ${\rm Div}\,A$ is locally a bounded measure. Then 
$$(\det A)^{\frac1{d-1}}\in L^1_{\rm loc}(\Omega).$$
\end{thm}

It is interesting to compare this statement with what we obtain when applying S. M\"uller's  Theorem in \cite{Mul} (see also Coifman \& all. \cite{CLMS}) to a vector field $u=\nabla\theta$ and $A:=\widehat{{\rm D}^2\theta}$. Theorem 1 in \cite{CLMS} tells us that if $\theta\in W^{2,d}_{\rm loc}$, then $f:=\det{\rm D}^2\theta$ belongs locally to the Hardy space ${\mathcal H}^1$ (a strict subspace of the $L^1_{\rm loc}$ that the H\"older inequality would give us). If moreover $\theta$ is convex, then $f\ge0$ and this amounts to saying that $f\log(1+ f)\in L^1_{\rm loc}$, which is M\"uller's statement. If instead we assume that every minor of ${\rm D}^2\theta$ of size $d-1$ is locally integrable (this is achieved for instance if $\theta\in W^{2,d-1}_{\rm loc}$), then $A\in L^1_{\rm loc}$ and Theorem \ref{th:gain} tells us that $(\det A)^{\frac1{d-1}}=f\in L^1_{\rm loc}$. Our result is actually implicit in  \cite{Mul}, where the inequality (2) meets our Theorem \ref{th:conv} (\ref{detconv}) when the vector field is a gradient, except for a non-optimal constant~; that inequality is attributed to H. Federer, Thm 4.5.9 (31) \cite{Fed}.

It would be interesting to understand the gain of integrability when $\theta\in W^{2,p}_{\rm loc}$ where $p\in(d-1,d)$.

\subsubsection{Application to the isoperimetric inequality}

Taking $A(x)\equiv I_d$, which is obviously a DPT, (\ref{detconv}) yields 
$$|\Omega|\le\frac1{d|S^{d-1}|^{\frac1{d-1}}}\,|\partial\Omega|^{\frac d{d-1}},$$
that is
$$\frac{|\Omega|}{|B_1|}\le\left(\frac{|\partial\Omega|}{|S_1|}\right)^{\frac d{d-1}}\,.$$

Although the proof above works only for convex domains, it can be adapted to general domains $E$. The following argument is due to G. de Philippis (personal communication). Choose a ball $\Omega$, which strictly contains $E$. Apply (\ref{convmeas}) to the tensor $A:={\bf1}_EI_d$, noticing the identity ${\rm Div}\,A=\nabla{\bf1}_E$. We obtain
$$|E|\le\frac1{d|S^{d-1}|^{\frac1{d-1}}}\,\|{\rm Div}\,A\|_{{\cal M}(\Omega)}^{\frac d{d-1}}=\frac1{d|S^{d-1}|^{\frac1{d-1}}}\,{\rm per}(E)^{\frac d{d-1}},$$
where ${\rm per}(E)$ is the perimeter of $E$ in the sense of Caccioppoli. This inequality is the isoperimetric one.



\bigskip

We shall see that the proof of Theorems \ref{th:per} and \ref{th:conv} are based upon mass transportation. A link between isoperimetric inequalities and mass transportation had already been noted by M. Gromov \cite{Gro}. However, Gromov's proof involves Knothe's map, whereas ours uses Brenier's map of optimal transport~; it is therefore closer to that of Figalli \& al. \cite{FMP}.

\section{Applications to gas dynamics}

We intend to apply or adapt Theorem \ref{th:conv} in a situation where the first independent variable is a time variable, and the other ones represent spatial coordinates. We therefore set $d=1+n$ and $x=(t,y)$ where $t\in\R$ and $y\in\R^n$. We  write a DPT blockwise
$$A(t,y)=\begin{pmatrix} \rho & m^T \\ m & S \end{pmatrix},$$
where $\rho\ge0$ and $m$ can be interpreted as the densities of mass and  linear momentum. We begin with an abstract result.
\begin{thm}\label{th:absfl}
Let $A$ be a DPT over a slab $(0,T)\times \R^n$. We assume 
$$A\in L^1((0,T)\times\R^n)\cap L^{\frac d{d-1}}_{loc}((0,T)\times\R^n).$$
There exists a constant $c_n$, depending only upon the space dimension (but neither on $T$, nor on $A$) such that, with the notations above
$$\int_0^Tdt\int_{\R^n}(\det A)^{\frac1n}dy\le c_n\left(\|m(0,\cdot)\|_{L^1(\R^n)}+\|m(T,\cdot)\|_{L^1(\R^n)}\right)\left(\int_{\R^n}\rho(0,y)\,dy\right)^{\frac1n}.$$
\end{thm}

\subsection{Euler equations} 

For a compressible, inviscid gas, the flux of momentum is given by $$S=\frac{m\otimes m}\rho+p I_n,$$
where the pressure $p\ge0$ is given by an equation of state. The latter is expressed in terms of the density $\rho$ (if the gas is barotropic or isentropic) or of the density and the temperature $\vartheta$ (adiabatic gas). In both cases, the Euler system ${\rm Div}\,A=0$ accounts for the conservation of mass and momentum, and is supplemented by an energy balance law
$$\partial_tE+{\rm div}_y\left[(E+p)\frac m\rho\right]\le0,\qquad E:=\frac{|m|^2}{2\rho}+\rho e,$$
where $e\ge0$ is the internal energy per unit mass. This inequality is an equality in the adiabatic case. Its main role is to provide an {\em a priori} energy estimate
$$\sup_{t\ge0}\int_{\R^n}E(t,y)\,dy\le E_0:=\int_{\R^n}E(0,y)\,dy,$$
whenever the total energy $E_0$ at initial time is finite. 

For reasonable equations of state, like those of a polytropic gas ($p=a\rho^\gamma$ for a constant $\gamma>1$) or a perfect gas ($p=(\gamma-1)\rho e$), the internal energy per unit volume $\rho e$ dominates the pressure: 
\begin{equation}\label{pCe}
p\le C\rho e
\end{equation} 
for some finite constant $C$. 

For a flow whose mass and energy are locally finite (a rather reasonable assumption), the tensor $A$ is locally integrable. Applying Theorem \ref{th:gain}, we deduce that $\rho^{\frac1n}p$ is locally integrable in space and time.  This is already a different and somehow better integrability than the ones obtained directly from the conservation of mass and energy decay.

If in addition the total energy $E_0$ is finite, we have that $S\in L^1((0,T)\times\R^n)$. If the total mass
$$M_0:=\int_{\R^n}\rho(0,y)\,dy$$
is finite too, then $A\in L^1((0,T)\times\R^n)$ (remark that the total mass remains constant in time).
Applying Theorem \ref{th:absfl} to the Euler system, we infer the estimate
$$\int_0^Tdt\int_{\R^n}\rho^{\frac1n}p\,dy\le 2c_n\,M_0^{\frac12+\frac1n}(2E_0)^{1/2}.$$
This inequality can be sharpened after remarking that the left-hand side does not depend upon the Galilean frame, while the right-hand side, more precisely $E_0$, does. We may replace in the inequality above the initial velocity $u_0=\frac m\rho(0,\cdot)$ by $u_0-\vec c$ where $\vec c$ is an arbitrary constant (this constant represents the velocity of a Galilean frame with respect to a reference frame). Eventually, we may choose the vector $\vec c$ which minimizes the resulting quantity
$$\int_{\R^n}\left(\frac12\,\rho_0|u_0-\vec c|^2+\rho_0 e_0\right)\,dx.$$
This yields the following result.
\begin{thm}\label{th:Euler} We assume that the equation of state implies (\ref{pCe}).

Consider an admissible (in the sense above) flow, solution of the Euler equations of a compressible fluid in $(0,T)\times \R^n$. We assume a finite mass $M_0$ and energy $E_0$, and that the tensor $A$ belongs to $L^{\frac d{d-1}}_{loc}((0,T)\times\R^n)$. Then the following estimate holds true:
\begin{equation}\label{rhopn}
\int_0^Tdt\int_{\R^n}\rho^{\frac1n}p\,dy\le 2c_n\,M_0^{\frac1n}D_0^{\frac12}, \end{equation}
where
$$D_0:=\int_{\R^n}\rho_0dy\int_{\R^n}(\rho_0|u_0|^2+2\rho_0e_0)\,dy-\left|\int_{\R^n}\rho_0u_0dy\right|^2.$$
\item 
\end{thm}

\paragraph{Remarks.}
\begin{itemize}
\item A careful examination gives the following value of the constant in (\ref{rhopn}):
$$c_n=\frac{(n+1)^{\frac1{2n}-\frac12}}{|S^n|^{\frac1n}\sqrt n}\,.$$
\item For full gas dynamics, the quantity  $D_0$ is an invariant of the flow. For a barotropic flow, the energy may decay, but the mass and linear momentum are preserved~; the corresponding quantity $D(t)$ is therefore non-increasing.
\item The kinetic part in $D_0$ can be written in a more symmetric way:
$$\int_{\R^n}\rho_0dy\int_{\R^n}\rho_0|u_0|^2dy-\left|\int_{\R^n}\rho_0u_0dy\right|^2=\frac12\int_{\R^n}\!\!\int_{\R^n}\rho_0(y)\rho_0(y')|u_0(y')-u_0(y)|^2dy\,dy',$$
in which the independence upon the choice of the Galilean frame becomes obvious. 
\item We did not make any local hypothesis about the equation of state. We did not even ask for hyperbolicity. Thus (\ref{rhopn}) could be used to control the blow-up for models with phase transition (Van der Waals gas). Our assumption (\ref{pCe}) is merely of a global nature. For instance, if the gas is barotropic, then $\rho\mapsto p,e$ are linked by $p=\rho^2e'$ and our assumption is just that
$$\rho\,\frac{de}{d\rho}\le C e$$
for some finite constant $C$.
\item When the solution is globally defined, we even have
\begin{equation}\label{globDel}
\int_0^\infty dt\int_{\R^n}\rho^{\frac1n}p\,dy\le 2c_n\,M_0^{\frac1n}D_0^{\frac12}.
\end{equation}
\item Our estimate shows that the fluid cannot  concentrate, unless $\rho^{\frac1n}p=O(\rho)$ as $\rho\rightarrow+\infty$. This rules out the so-called {\em delta-shocks} for most of the reasonable equations of state.
\end{itemize}

\paragraph{Polytropic gas.}
When $p(\rho)={\rm cst}\cdot\rho^\gamma$ with adiabatic constant $\gamma>1$, (\ref{rhopn}) is an estimate of $\rho$ in $L^{\gamma+\frac1n}_{t,y}$, which up to our knowledge is new. 
Combining this with the estimates of $\rho$ in $L^\infty_t(L^1_y)$ (conservation of total mass) and in $L^\infty_t(L^\gamma_y)$ (decay of total energy), and using the H\"older inequality, we infer estimates of $\rho$ in $L^q_t(L^r_y)$ for every $(q,r)$ such that the point $\left(\frac1q,\frac1r\right)$ lies within the triangle whose vertices are
$$(0,1),\quad\left(0,\frac1\gamma\right)\quad\hbox{and}\quad\left(\frac n{n\gamma+1}\,,\frac n{n\gamma+1}\right).$$
A similar interpolation argument, which involves the decay of energy, ensures that 
$$\int_0^T\left(\int_{\R^n}\rho^\alpha|m|\,dy\right)^2dt<\infty,\qquad\alpha:=\frac12\left(\frac1n+\gamma-1\right).$$
When $T=+\infty$, (\ref{globDel}) can be compared with other dispersion estimates, for instance (see \cite{Ser_eter})
$$\int_{\R^n}\rho^\gamma dy=O\left((1+t)^{-n(\gamma-1)}\right),$$
when the gas has finite inertia
$$I_0:=\int_{\R^n}\rho(0,y)\,\frac{|y|^2}2\,dy.$$

\paragraph{Perfect non-isentropic gas.} When $p=(\gamma-1)\rho e$, a similar argument yields an estimate of $\rho^{1+\frac1{nq}}e^r$ in $L^q_t(L^1_x)$, whenever $1\le q\le\infty$ and $r-1\le\frac1q\le r$.

\paragraph{Euler--Fourier system.}

The Euler--Fourier system governs the motion of an inviscid but heat-conducting gas. The only difference with the Euler system is that the conservation law of energy incorporates a dissipative term ${\rm div}_y(\kappa\nabla_y\vartheta)$, where $\vartheta$ is the temperature and $\kappa>0$ the thermal conductivity. Because the conservation of mass and momentum still writes ${\rm Div}\,A$ with the same $A$ as before, and because the total energy is conserved, Theorem \ref{th:Euler} applies to this case.

On the contrary, our theorem does not apply to the Navier--Stokes system for a compressible fluid, because then the divergence-free tensor is not positive semi-definite.

\paragraph{The role of Estimate (\protect\ref{rhopn}).}
Theorem \ref{th:Euler}  is an {\em a priori} estimate which suggests a functional space where to search for admissible solutions of the Euler equation. For finite initial mass and energy, one should look for a flow satisfying
the following three requirements: -- the total mass is conserved, -- the total energy is a non-increasing function of time (a constant in the adiabatic case), -- and $\rho^{\frac1n}p\in L^1_{t,y}$. 

To this end, the construction of a solution to the Cauchy problem should involve an approximation process which is consistent with these estimates. For this purpose, a vanishing viscosity approach (say, the compressible Navier--Stokes equation) does not seem suitable. As we shall see below, the Boltzmann equation is more appropriate, but this observation just shifts the consistency problem to an other level. An other approach is to design numerical schemes, which are consistent with the Euler equations and meanwhile with the above requirements. There exist several schemes that preserve the symmetric positive structure, for instance the Lax--Friedrichs and Godunov schemes in space dimension one, or their muti-dimensional variants. However they provide approximations for which the mass of the Radon measure ${\rm Div}\, A^{\Delta t,\Delta y}$ tends to $+\infty$ as $\Delta t,\Delta y\rightarrow0$. The second part  of Theorem \ref{th:conv} yields
$$\int_{\Omega}(\det A^{\Delta t,\Delta y})^{\frac1{d-1}}dy\,dt\le\frac1{d|S^{d-1}|^{\frac1{d-1}}}\,\left(\|A^{\Delta t,\Delta y}\vec n\|_{L^1(\partial\Omega)}+\|{\rm Div}\,A^{\Delta t,\Delta y}\|_{{\cal M}(\Omega)}\right)^{\frac d{d-1}},$$
where the right-hand side tends to $+\infty$ when $\Delta t,\Delta y\rightarrow0$. Thus it is unclear whether the limit of such schemes satisfies the estimate (\ref{rhopn}).

Notice that we must not require the integrability $A\in L^{\frac d{d-1}}$, which is only a technical need for our proof. As a matter of fact, the various entries $a_{ij}$ have distinct physical dimensions, so that such an integrability hardly makes sense. On the contrary, $\det A$ is a well-defined quantity from the physical point of view. 

We also point out that, although our new estimate is a genuine improvement, it is still not sufficient to ensure the local integrability of the energy flux
$$\left(\frac12\rho|u|^2+\rho e+p\right)u,$$
and therefore to give sense to the conservation law of energy.

\subsection{Self-similar flows}

We now consider the problem (\ref{RP}) in space dimension $n$. The tensor $A=\rho v\otimes v+pI_n$ (recall that $v$ is the pseudo-velocity), though positive semi-definite, is not a DPT. The second part of Theorem \ref{th:conv}, plus the formula $\det A=p^{n-1}(p+\rho|v|^2)$, yield
\begin{eqnarray}
\label{inRP}
\int_\Omega p(p+\rho|v|^2)^{\frac1{n-1}}d\xi & \le & \frac1{n|S^{n-1}|^{\frac1{n-1}}}\,\left(\|p\vec n+\rho(v\cdot\vec n)v\|_{L^1(\partial\Omega)}+(n+1)\int_\Omega\rho |v|\,d\xi\right)^{\frac{n}{n-1}} \\
\nonumber & \le & \frac1{n|S^{n-1}|^{\frac1{n-1}}}\,\left(\int_{\partial\Omega}(p+\rho|v|^2)\,ds(\xi)+(n+1)\int_\Omega\rho |v|\,d\xi\right)^{\frac{n}{n-1}}
\end{eqnarray}
for every convex subdomain $\Omega$. For a ball $B_r$ of radius $r$ and arbitrary center, this writes
\begin{equation}\label{inRPball}
\Xint-_{B_r} p(p+\rho|v|^2)^{\frac1{n-1}}d\xi \le \left(\Xint-_{\partial B_r}(p+\rho|v|^2)\,ds(\xi)+\frac{n+1}n\,r\Xint-_{B_r}\rho |v|\,d\xi\right)^{\frac{n}{n-1}}
\end{equation}
Remark that, contrary to the situation of the Cauchy problem, we do not have the freedom to choose among equivalent coordinate frames. There is no improvement of (\ref{inRP}) or (\ref{inRPball}) similar  to (\ref{rhopn}).

\paragraph{Riemann problem.}
The Riemann problem is a special form of the Cauchy problem, where the initial data (density, momentum, energy) is positively homogeneous of degree zero~; for instance, the initial density has the form $\rho_0(\frac y{|y|})$\,. In practice, we suppose that the physical space $\R^n$ is partitioned into conical cells with polygonal sections, and that the data is constant in each cell. Such a data depends on finitely many parameters.

Because the Euler equations are PDEs of homogeneous degree one, the admissible solution, whether there exists a unique one, must be self-similar too. The density satisfies $\rho(t,x)=\bar\rho(\frac xt)$\, and so on. Denoting $\xi=\frac xt$ the self-similar variable, every conservation law $\partial_tf+{\rm div}_yq=0$ becomes ${\rm div}_\xi q=\xi\cdot\nabla_\xi f$. For instance, droping the bars, we have ${\rm div}_\xi(\rho u)=\xi\cdot\nabla \rho$. These new equations involve explicitly the independent variable $\xi$, but the introduction of the {\em pseudo-velocity} $v(\xi):=u(\xi)-\xi$ allows us to get rid of it. In terms of $\rho, v,p,e$ and $\xi$-derivatives, the reduced Euler system becomes
\begin{equation}\label{clRP}
{\rm div}(\rho v)+n\rho=0,\qquad{\rm Div}(\rho v\otimes v)+\nabla p+(n+1)\rho v=0
\end{equation}
and
\begin{equation}\label{enerRP}
{\rm div}\left((\frac12\rho|v|^2+\rho e+p)v\right)+\left(\frac n2+1\right)\rho |v|^2+n(\rho e+p)=0.
\end{equation}
The initial data to the Riemann problem becomes a data at infinity for the reduced system. Let us mention that for an isentropic flow, (\ref{enerRP}) is not an equation but merely an inequation, which plays the role of an entropy inequality.

\bigskip

The $3$-dimensional RP is still widely open. We therefore limit ourselves to the $2$-dimensional case ($n=2$). The tools and strategy for the analysis of the Riemann problem are described in the review paper \cite{Ser_refl}. The plane splits into a compact subsonic region $\Omega_{\rm sub}$ and its complement the supersonic domain $\Omega_{\rm sup}$. {\em Subsonic} means that $|v|\le{\bf c}$ where $\bf c$ is the sound speed, a function of the internal variables $\rho$ and $p$. In the supersonic region, the system is of hyperbolic type and one can solve a kind of Cauchy problem, starting from the data at infinity. This Cauchy problem has an explicit solution outside some ball $B_R(0)$. It is made of constant states separated by simple waves depending only on one coordinate~; these waves are shocks, rarefaction waves and/or contact discontinuities. An {\em a priori} estimate of the radius $R$ is available. The situation in the rest of the supersonic region may be more involved, with genuinely $2$-D interactions of simple waves~; even the interface between $\Omega_{\rm sup}$ and $\Omega_{\rm sub}$ is not fully explicit, a part of it being a free boundary. But these facts do not raise obstacles for the following calculations.

The conservation laws of mass and energy allow us to establish two {\em a priori} estimates. On the one hand, we have (recall that $n=2$)
$$2\int_{B_R(0)}\rho\,d\xi=-\int_{B_R(0)}{\rm div}(\rho v)\,d\xi=-\int_{S_R(0)}\rho v\cdot\vec n\,d\xi,$$
where the last integral is computed explicitly because of our knowledge of the solution over $S_R$. On the other hand, we have
$$2\int_{B_R(0)}(\rho |v|^2+\rho e+p)\,d\xi\le-\int_{S_R(0)}(\frac12\rho|v|^2+\rho e+p) v\cdot\vec n\,d\xi,$$
where again the right-hand side is known explicitly. In the non-isentropic case, we also have a minimum principle for the physical entropy $s$, which is nothing but the second principle of thermodynamics: $s\ge s_{\min}$ where $s_{\min}$ is the minimum value of $s$ in the data. Let us assume a polytropic gas ($p={\rm cst}\cdot\rho^\gamma$) or a perfect gas ($p=(\gamma-1)\rho e$) law. In the latter case, we have $p\ge(\gamma-1)e^{s_{\min}}\rho^\gamma$. Therefore the energy estimate yields an upper bound for
\begin{equation}\label{firstestRP}
\int_{B_R(0)}\rho^\gamma d\xi\qquad\hbox{and}\qquad\int_{B_R(0)}\rho|v|^2d\xi.
\end{equation}
In particular, a so-called {\em Delta-shock} cannot take place in this situation.

These estimates can be completed by applying  (\ref{inRP}) to the tensor $A=\rho v\otimes v+pI_2$ in the ball $B_R(0)$. To this end, we show that the right-hand side is fully controled. On the one hand, the boundary integral
$$\int_{S_R(0)}(p+\rho|v|^2)\,ds(\xi)$$
is estimated explicitly as before.  On the other hand, the last integral is bounded by
$$\left(\int_{B_R(0)}\rho\,d\xi\right)^{\frac12}\left(\int_{B_R(0)}\rho|v|^2\,d\xi\right)^{\frac12},$$
where both factors have been estimated previously. We therefore obtain an estimate of
\begin{equation}\label{imprRP}
\int_{B_R(0)}\rho^{2\gamma}d\xi\qquad\hbox{and}\qquad\int_{B_R(0)}\rho^{\gamma+1}|v|^2d\xi.\end{equation}
This integrability is significantly better than that in (\ref{firstestRP}).

\subsection{Relativistic gas dynamics}

In the Minkowski space-time $\R^{1+n}$ of special relativity, an isentropic gas is governed by the Euler system (see \cite{MU})
\begin{eqnarray*}
\partial_t\left(\frac{\rho c^2+p}{c^2-|v|^2}\,-\,\frac p{c^2}\right)+{\rm div}_y\left(\frac{\rho c^2+p}{c^2-|v|^2}\,v\right) & = & 0, \\
\partial_t\left(\frac{\rho c^2+p}{c^2-|v|^2}\,v\right)+{\rm Div}_y\left(\frac{\rho c^2+p}{c^2-|v|^2}\,v\otimes v\right)+\nabla p & = & 0,
\end{eqnarray*}
where the constant $c>0$ is the speed of light. Here $\rho$ is the mass density at rest, and $p$ is the pressure. The fluid velocity is constrained by $|v|<c$.

It is a simple exercise to verify that the stress-energy tensor 
$$A=\begin{pmatrix} \frac{\rho c^2+p}{c^2-|v|^2}\,-\,\frac p{c^2} & \frac{\rho c^2+p}{c^2-|v|^2}\,v^T  \\
\frac{\rho c^2+p}{c^2-|v|^2}\,v & \frac{\rho c^2+p}{c^2-|v|^2}\,v\otimes v + pI_n
\end{pmatrix}$$
is positive semi-definite. Our Theorems \ref{th:conv} and \ref{th:absfl} therefore apply. What is perhaps surprising is that the determinant of $A$ is unchanged~! Its value is still $\rho p^n$. We 
infer
$$\int_0^Tdt\int_{\R^n}\rho^{\frac1n}p\,dy\le c_n\left(\left.\left.\int_{\R^n}\frac{\rho c^2+p}{c^2-|v|^2}\,|v|\,dy\right|_{t=0}+(\hbox{same)}\right|_{t=T}\right)\left(\int_{\R^n}\frac{\rho c^4+p|v|^2}{c^2(c^2-|v|^2)}\,dy\right)^{\frac1n}_{t=0}.$$
We warn the reader that mass and energy are related to each other in relativity theory. The last integral in the inequality above accounts for both. We denote below its value $\mu_0$.

\bigskip

Suppose an equation of state of the form $p=a^2\rho$, where $a>0$ is a constant. When $a^2=\frac{c^2}3$\,, this follows directly from the Stefan--Boltzmann law for a gas in thermodynamical equilibrium, as discussed page 12 of A. M. Anile's book \cite{Ani}. Then the contribution of the momentum can be estimated after using $|v|\le\frac1{2ca}(c^2+a^2|v|^2)$~:
$$\int_{\R^n}\frac{\rho c^2+p}{c^2-|v|^2}\,|v|\,dy\le\frac{c^2+a^2}{2a}\,\mu_0.$$
We deduce the {\em a priori} estimate
\begin{equation}\label{estRR}
\int_0^Tdt\int_{\R^n}\rho^{1+\frac1n}\,dy\le c_n\frac{c^2+a^2}{a^3}\,\mu_0^{1+\frac1n}.
\end{equation}

\subsection{Kinetic equations}

We now turn to the class of kinetic equations
\begin{equation}\label{kineq}
(\partial_t+v\cdot\nabla_y)f(t,y,v)=Q[f(t,y,\cdot)]
\end{equation}
where $Q$ is compatible with the minimum principle $f\ge0$ and with the conservation of mass, momentum and energy. This includes the Boltzman equation, the BGK model and most of the discrete velocity models. Then we apply Theorem \ref{th:conv} to the non-negative tensor
$$A(t,y):=\int_{\R^n}f(t,y,v)\begin{pmatrix} 1 \\ v \end{pmatrix}\otimes\begin{pmatrix} 1 \\ v \end{pmatrix}\,dv.$$
The following result is a far-reaching extension of an estimate that J.-M. Bony \cite{Bony} obtained for a one-dimensional discrete velocity model.
\begin{thm}\label{th:kin}
Consider an admissible flow of a kinetic equation of the form (\ref{kineq}). Assume a finite mass an energy
$$M_0=\int_{\R^n}dy\int_{\R^n}f_0(y,v)\,dv,\qquad E_0=\int_{\R^n}dy\int_{\R^n}f_0(y,v)\frac{|v|^2}2\,dv,$$
and that the moments 
$$\rho(t,y)=\int_{\R^n}f(t,y,v)\,dv,\qquad \Tr S(t,y)=\int_{\R^n}f(t,y,v)|v|^2\,dv$$
belong to $L^{\frac d{d-1}}_{loc}((0,T)\times\R^n)$.
Then the following estimate holds true:
\begin{equation}\label{desfois}  \int_0^Tdt\int_{\R^n}dy\left(\frac1{d!}\int_{\R^n}^{\otimes(n+1)}f(t,y,v^0)\cdots f(t,y,v^n)(\Delta(v^0,\ldots,v^n))^2dv^0\cdots dv^n\right)^{\frac1n}\le 2c_n\,M_0^{\frac1n}D_0^{1/2}, 
\end{equation}
where
$$\Delta(v_0,\ldots,v_n):=\begin{vmatrix} 1 & \cdots & 1 \\ v^0 & \cdots & v^n \end{vmatrix}$$
is $n!$ times the volume of the simplex spanned by $(v^0,\ldots,v^n)$, and
$$D_0=\frac12\int_{\R^n}^{\otimes 4}f_0(y,v)f_0(y',v')|v'-v|^2dydy'dvdv'.$$
\end{thm}

Again, this estimate suggests to narrow the functional space where to search for a solution. Besides the usual constraints
$$\sup_t\int_{\R^n}\!\!\int_{\R^n}(1+|x|^2+|v|^2+\log^+ f)\,f\,dv\,dx<\infty,$$
we should impose that the expression
$$I_T:=\int_0^Tdt\int_{\R^n}dy\left(\frac1{d!}\int\!\cdots\!\int_{\R^n}^{\otimes(n+1)}f(t,y,v^0)\cdots f(t,y,v^n)(\Delta(v^0,\ldots,v^n))^2dv^0\cdots dv^n\right)^{\frac1n}$$
be finite. An open problem is to understand the physical meaning of $I_T$.
\paragraph{Comments.}
\begin{itemize}
\item The $d\times d$ determinant $\Delta(v_0,\ldots,v_n)$ vanishes precisely  when the points $v^0,\ldots,v^n$ are affinely dependent in the space $\R^n$, therefore are non generic. The estimate (\ref{desfois}) tends to force the support of $f(t,y,\cdot)$ to keep close to some affine hyperplane $\Pi(t,y)$.
\item Of course, Boltzman's $H$-theorem, which tends to force $f(t,y,\cdot)$ to be close to a Maxwellian distribution, has the opposite effect. The combination of both estimates is expected to produce a nice control of the density $f$.
\item Our estimate controls the $(t,y)$-integrability of an expression homogeneous in $f$ of degree $1+\frac1n$. This is slightly but strictly better than the controls given by the mass and energy (both linear in $f$) or by the H-Theorem (control in $f\log f$). The price to pay is an integration in the time variable ; this looks like what happens in Strichartz estimates for dispersive PDEs.

The little gain in integrability raises the question whether the Boltzmann equation admits weak solutions for large data, and not only renormalized ones. Using this gain, C. Cercignani \cite{Cer} proved the existence of weak solutions to the Cauchy problem in dimension $n=1$.
\item If we had just applied the Jensen inequality, the exponent in (\ref{desfois}) would have been $\frac1{n+1}\,$, and the $(t,y)$-integrand should be homogeneous of degree $1$, conveying an information already contained in the mass and energy.
\end{itemize}

\subsubsection{Renormalized solutions}

The existence of distributional solutions to the Cauchy problem for the Boltzmann equation has not yet been proved, except in space dimension $n=1$. R. DiPerna \& P.-L. Lions \cite{DPL} proved the existence of a weaker notion of solutions, called {\em renormalized}. We shall not even give a precise definition of this notion, but we content ourselves to recall that it implies the conservation of mass and a weak form of the conservation of momentum, in the sense that 
\begin{equation}\label{Euldef}
\partial_t\int_{\R^n}f\,dv+{\rm div}_y\int_{\R^n}fv\,dv=0,\qquad\partial_t\int_{\R^n}fv\,dv+{\rm Div}_y\left(\int_{\R^n}fv\otimes v\,dv+\Sigma\right)=0.
\end{equation}
Compared to what is formally expected, the second equation above contains an additional term, called the {\em defect measure} $\Sigma$, which takes values in ${\bf Sym}_n^+$~; see \cite{LM}. Finally, it is known that the total energy
$$E(t):=\frac12\int_{\R^n}\!\!\int_{\R^n}f(t,y,v)|v|^2dv\,dy+\frac12\,\int \Tr\Sigma(t,\cdot)$$
is a non-increasing function of time and satisfies
$$E(t)\le E_0:=\frac12\int_{\R^n}\!\!\int_{\R^n}f_0(y,v)|v|^2dv\,dy.$$

The equations (\ref{Euldef}) can be recast by saying that the following tensor is a DPT
$$A=\begin{pmatrix} \int_{\R^n}f\,dv & \int_{\R^n}fv^T\,dv \\ \int_{\R^n}fv\,dv & \int_{\R^n}fv\otimes v\,dv+\Sigma \end{pmatrix}.$$
The components $\rho,m$ of $A\vec e_t$ are still the mass density and linear momentum. Theorem \ref{th:absfl} yields an inequality
$$\int_0^Tdt\int_{\R^n}(\det A)^{\frac1n}dy\le c_nM_0^{\frac1n}(\|m(0)\|_{L^1(\R^n)}+\|m(T)\|_{L^1(\R^n)}),$$
from which we may extract two informations, using the monotonicity of the determinant. On the one hand, we have
$$\begin{pmatrix} \int_{\R^n}f\,dv & \int_{\R^n}fv^T\,dv \\ \int_{\R^n}fv\,dv & \int_{\R^n}fv\otimes v\,dv \end{pmatrix}\le A,$$
from which we infer the same estimate (\ref{desfois}) as in the case of distributional solutions. On the other hand, the Schur complement formula  (see \cite{Ser_mat} Proposition 3.9) gives
$$\det A=\rho\det\left(\int_{\R^n}fv\otimes v\,dv-\frac1{\int_{\R^n}f\,dv}\int_{\R^n}fv\,dv\otimes\int_{\R^n}fv\,dv+\Sigma\right)\ge\rho\det\Sigma,$$
because the tensor
$$\int_{\R^n}f\,dv\int_{\R^n}fv\otimes v\,dv-\int_{\R^n}fv\,dv\otimes\int_{\R^n}fv\,dv$$
is positive semi-definite.
We infer an estimate of the defect measure against the mass density:
\begin{equation}\label{intsig}
\int_0^Tdt\int_{\R^n}(\rho\det \Sigma)^{\frac1n}\le c_n'M_0^{\frac1n}D_0^{\frac12}.
\end{equation}

Notice that, because $\Sigma(t,\cdot)$ is a Radon measure taking values in ${\bf Sym}_n^+$ and $\det^{\frac1n}$ is homogeneous of degree one over this cone, the expression $(\det \Sigma(t,\cdot))^{\frac1n}$ makes sense as a bounded measure.

\section{Convex domain}\label{s:convex}

\subsection{Proof of Theorem \protect\ref{th:conv}}

In this paragraph, we consider a DPT over a bounded convex domain $\Omega$. To prove Theorem \ref{th:conv}, it is enough to consider the case where $A$ is uniformly positive definite: just replace $A(x)$ by $A(x)+\delta I_d$ with $\delta>0$ (such a tensor is still a DPT) and then pass to the limit as $\delta\rightarrow0+$. 

From now on, we therefore assume that $\Omega$ has a smooth boundary and that $A(x)\ge\delta I_d$ for some $\delta>0$ independent of $x$.

\bigskip

Let $f$ denote the function $(\det A)^{\frac1{d-1}}$. One has $f=(f\det A)^{\frac1d}$. The density of ${\cal C}^\infty(\overline\Omega)$ in $L^1(\Omega)$ provides
a sequence of smooth functions $f_\epsilon:\Omega\rightarrow\R$ that satisfies the following requirements. To begin with, $\frac\delta2\le f_\epsilon(x)\le C_\epsilon$ for every $x$, where $C_\epsilon$ is a finite constant depending on $\epsilon$. Then
$$\int_\Omega f_\epsilon(x)\,dx=\int_\Omega f(x)\,dx,$$
and finally
$$\|f_\epsilon-f\|_{L^1(\Omega)}\stackrel{\epsilon\rightarrow0}\longrightarrow0.$$
From the latter, we deduce that $f_\epsilon^{1/d}\rightarrow f^{1/d}$ in $L^d(\Omega)$ and therefore $f_\epsilon^{1/d}f^{1-1/d}\rightarrow f$ in $L^1(\Omega)$. It will thus be enough to estimate
$$\int_\Omega f_\epsilon^{1/d}f^{1-1/d}\,dx=\int_\Omega (f_\epsilon\det A)^{1/d}dx.$$
To do so, we consider the ball $B_r=B_r(0)$, centered at the origin, whose volume equals the integral of $f$ (that is, that of $f_\epsilon$) over $\Omega$. A theorem due to Y. Brenier (see Theorem 2.12 in \cite{Vil}, or Theorem 3.1 in \cite{dPF}) ensures the existence and uniqueness of an optimal transport from the measure $f_\epsilon(x)\,dx$ to the Lebesgue measure over $B_r$. This transport is given by a gradient map $\nabla \psi_\epsilon$, which is the solution of the Monge--Amp\`ere equation
$$\det\nabla^2\psi_\epsilon=f_\epsilon\qquad\hbox{in }\Omega$$
such that $\psi_\epsilon$ is convex and $\nabla\psi_\epsilon(\Omega)=B_r$~; see Theorem 4.10 of \cite{Vil} or Theorem 3.3 of \cite{dPF}. Finally, $\psi_\epsilon$ is a smooth function (Theorem 4.13 of \cite{Vil}). In particular, the image of the boundary $\partial\Omega$ under $\nabla\psi_\epsilon$ is the sphere $S_r$.

We therefore have 
$$(f_\epsilon\det A)^{\frac1d}=(\det A\cdot\det\nabla^2\psi_\epsilon)^{\frac1d}=(\det (A\nabla^2\psi_\epsilon))^{\frac1d}.$$
Let $\lambda_1(x),\ldots,\lambda_d(x)$ be the spectrum of the matrix $A\nabla^2\psi_\epsilon$. This is not a symmetric matrix, but because it is the product of a positive definite matrix and a positive semi-definite one, it is diagonalisable with non-negative real eigenvalues: the $\lambda_j$'s are real and $\ge0$ (Proposition 6.1 in \cite{Ser_mat}). Applying the geometric-arithmetic mean inequality (AGI), we have
$$(\det(A\nabla^2\psi_\epsilon))^{\frac1d}=\left(\prod_{j=1}^d\lambda_j\right)^{\frac1d}\le\frac1d\sum_{j=1}^d\lambda_j=\frac1d\,\Tr(A\nabla^2\psi_\epsilon).$$

Because $A$ is divergence-free, one has $\Tr(A\nabla^2\psi_\epsilon)={\rm div}(A\nabla\psi_\epsilon)$. We infer
$$\int_\Omega f_\epsilon(x)^{\frac1d}f(x)^{1-\frac1d}dx\le\frac1d\int_\Omega{\rm div}(A\nabla\psi_\epsilon)\,dx=\frac1d\int_{\partial\Omega}(A\nabla\psi_\epsilon)\cdot\vec n\,ds(x)=\frac1d\int_{\partial\Omega}(A\vec n)\cdot\nabla\psi_\epsilon\,ds(x).$$
Because $\nabla\psi_\epsilon$ takes values in $B_r$, there comes
$$\int_\Omega f(x)\,dx=\lim_{\epsilon\rightarrow0}\int_\Omega f_\epsilon(x)^{\frac1d}f(x)^{1-\frac1d}dx\le\frac rd\,\|A\vec n\|_{L^1(\partial\Omega)}.$$
We complete the proof of the theorem by the calculation of the radius $r$~:
$$\frac{r^d}d\,|S^{d-1}|=|B_r|=\int_\Omega f(x)\,dx.$$

\bigskip

If instead ${\rm Div}\,A$ is a bounded measure, then we have $\Tr(A\nabla^2\psi_\epsilon)={\rm div}(A\nabla\psi_\epsilon)-({\rm Div}\,A)\nabla\psi_\epsilon$. The same calculation yields the bound
$$\int_\Omega f(x)\,dx\le\frac rd\,\left(\|A\vec n\|_{L^1(\partial\Omega)}+\|{\rm Div}\,A\|_{{\cal M}(\Omega)}\right)$$
and the conclusion follows.

\paragraph{About the symmetry assumption:} in the proof above as well as in other forthcoming proofs, we use the property that for a positive semi-definite $d\times d$ matrix $A$, one has $(\det(AS))^{1/d}\le\frac1d\,\Tr(AS)$ for every $S\in{\bf SPD}_d$. For this to hold true, the assumption on $A$ is not only sufficient, it is also necessary: if $M\in{\bf M}_d(\R)$ with $\det M>0$ is such that $(\det(MS))^{1/d}\le\frac1d\,\Tr(MS)$ for every $S\in{\bf SPD}_d$, then $M\in{\bf SPD}_d$. To see this, just use the polar factorization $M=QH$ ($Q$ orthogonal and $H\in{\bf SPD}_d$) and choose $S=H^{-1}$. One obtains $1\le\frac1d\,\Tr Q$, which implies $Q=I_d$ because the $Q$ is diagonalizable with eigenvalues of unit modulus.

\subsection{The equality case: proof of Proposition \protect\ref{p:eqcase}}

Since we assume that $f$ is smooth and bounded below and above, we may take $f_\epsilon=f$. Let us examine the proof above. In order to have equalities everywhere, we need in particular that the AGI be an equality, that is the $\lambda_j$'s be equal to each other. Then the diagonalisable matrix $A\nabla^2\psi$, with only one eigenvalue $\lambda(x)$, must equal $\lambda(x)I_d$. In other words, there is a scalar field $a>0$ such that $A(x)=a(x)\widehat{\nabla^2\psi}$. In particular $\widehat{\nabla^2\psi}$ is positive definite. Because both $A$ and $\widehat{\nabla^2\psi}$ are divergence-free (Lemma \ref{l:widehat}), we find that $(\widehat{\nabla^2\psi})\nabla a=0$, that is $\nabla a=0$. Thus $a$ is a constant. Up to replacing $\psi$ by $a^{-1/(d-1)}\psi$, we infer that $A=\widehat{\nabla^2\psi}$. By construction the image of $\Omega$ by $\nabla\psi$ is a ball centered at the origin.

\bigskip

Conversely, if $\psi$ is such a convex function, and $A(x):=\widehat{\nabla^2\psi}$, then we know that $A$ is a DPT. Let us examine the calculations of the previous paragraph. There is no need of an $f_\epsilon$, we can just keep $f$ itself. Likewise, we take $\psi$ instead of $\psi_\epsilon$. Because $A\nabla^2\psi=(\det \nabla^2\psi)I_d$, the AGI is actually an equality and we have
$$\int_\Omega f(x)\,dx=\frac1d\int_\Omega{\rm div}(A\nabla\psi)\,dx=\frac1d\int_{\partial\Omega}(A\nabla\psi)\cdot\vec n\,dx=\frac1d\int_{\partial\Omega}(A\vec n)\cdot\nabla\psi\,dx.$$
We claim that $A\vec n$ and $\nabla\psi$ are positively colinear along the boundary. It amounts to proving that $\vec n$ and $A^{-1}\nabla\psi$ are so. But the latter vector equals 
$$\frac1{\det \nabla^2\psi}\,\nabla^2\psi\,\nabla\psi=\frac1{\det \nabla^2\psi}\,\nabla\left(|\nabla\psi|^2\right).$$
Because $|\nabla\psi|^2$ is $\le r$ everywhere, but equals $r$ on $\partial\Omega$, its gradient is normal to the boundary and points outward. This proves the claim.

We therefore have $(A\vec n)\cdot\nabla\psi=|A\vec n|\cdot|\nabla\psi|=r|A\vec n|$ over $\partial\Omega$, and we infer
$$\int_\Omega f(x)\,dx=\frac rd\,\|A\vec n\|_{1,\partial\Omega}.$$
This ends the proof of the proposition.

\subsection{Proof of Proposition \protect\ref{p:suppK}}

Because $A(x)$ is positive semi-definite, $\vec n^TA\vec n=0$ implies $A\vec n=0$. This ensures that the extension of $A$ to $\R^d$ by $A\equiv0$ over $\R^d\setminus\Omega$, is still a DPT over $\R^d$. Let us denote it $A$, which is compactly supported. Let $\phi_\epsilon$ be a non-negative mollifier and set $A^\epsilon=\phi_\epsilon*A$. This is a compactly supported DPT, of class $C^\infty$. Its Fourier transform is therefore in the Schwartz class. The divergence-free constraint translates into ${\cal F}A^\epsilon(\xi)\xi\equiv0$. Taking $\eta\in S^{d-1}$ and $r>0$, we have ${\cal F}A^\epsilon(r\eta)\eta=0$. Letting $r\rightarrow0+$, we obtain ${\cal F}A^\epsilon(0)\eta=0$. In other words ${\cal F}A^\epsilon(0)=0_d$, that is
$$\int_{\R^d}A^\epsilon(x)\,dx=0_d.$$
From there, the non-negativity of $A^\epsilon(x)$ for all $x$ implies $A^\epsilon\equiv0_d$. Passing to the limit as $\epsilon\rightarrow0+$, we infer $A\equiv0_d$.

\bigskip

\paragraph{Proof of Theorem \protect\ref{th:gain}.}
Let $B$ be a ball such that $\bar B\subset\Omega$.  As in the previous section, we may assume that $A$ is uniformly positive definite: $A(x)\ge\delta I_d$ for almost every $x$. We begin by mollifying $A$, defining $A_\epsilon=\phi_\epsilon*A$, where 
$$\phi_\epsilon(x)=\frac1{\epsilon^d}\,\phi\left(\frac x\epsilon\right),\qquad\phi\in{\cal D}^+(\R^d)\quad\hbox{and}\quad\int_{\R^d}\phi(x)\,dx=1.$$
This makes sense in $B$ whenever $\epsilon>0$ is small enough that $B+B_\epsilon\subset\Omega$. The resulting $A_\epsilon$ is a smooth, uniformly positive symmetric tensor in $B$.

Let $\chi\in{\mathcal D}^+(B)$ be given. We apply (\ref{convmeas}) to the non-negative tensor $\chi A_\epsilon$ over the domain $B$. Because $\chi$ is compactly supported, we obtain
\begin{equation}\label{chiAeps}
\int_B\chi^{\frac d{d-1}}(\det A_\epsilon)^{\frac1{d-1}}dx\le\frac1{|S^{d-1}|^{\frac1{d-1}}}\|{\rm Div}(\chi A_\epsilon)\|_{{\mathcal M}(B)}^{\frac d{d-1}}.
\end{equation}
Because ${\rm Div}(\chi A_\epsilon)=\chi {\rm Div}(A_\epsilon)+A^\epsilon\nabla\chi$ and 
$$\|{\rm Div}\,A_\epsilon\|_{{\mathcal M}(B)}\le\|{\rm Div}\, A\|_{{\mathcal M}(B+B_\epsilon)},\qquad \|A_\epsilon\|_{L^1(B)}\le\|A\|_{L^1(B+B_\epsilon)},$$
the right-hand side of (\ref{chiAeps}) remains bounded as $\epsilon\rightarrow0+$. Because $A_\epsilon\rightarrow A$ in $L^1(B)$, we have, up to the extraction of a sub-sequence, $A_\epsilon(x)\rightarrow A(x)$ almost everywhere and therefore $\det A_\epsilon(x)\rightarrow\det A(x)$. Passing to the limit in (\ref{chiAeps}) and using Fatou's Lemma, we obtain that
$$\int_B\chi^{\frac d{d-1}}(\det A)^{\frac1{d-1}}dx<\infty.$$
This proves the theorem.

\section{Periodic tensors: proofs}\label{s:per}

\subsection{Proof of Theorem \protect\ref{th:per}}

\paragraph{Reduction.}
As in the previous section, we may assume that $A$ is uniformly positive definite: $A(x)\ge\delta I_d$ for almost every $x$. Also, we may approximate $A$ by a smooth DPT $A_\epsilon=\phi_\epsilon*A$ as above.
This $A_\epsilon$ is smooth and still satisfies $A_\epsilon(x)\ge\delta I_d$. It converges towards $A$ in $L^1_{loc}$, and we may assume that $A_\epsilon(x)\rightarrow A(x)$ almost everywhere. In particular, $(\det A_\epsilon)^{\frac1{d-1}}$ converges almost everywhere towards $(\det A)^{\frac1{d-1}}$. If we know that $(\det A_\epsilon)^{\frac1{d-1}}\in L^1(\R^d/\Gamma)$ and that 
the inequality (\ref{Jdetper}) holds true for $A_\epsilon$, then we may pass to the limit and Fatou's Lemma implies that  (\ref{Jdetper}) holds true for $A$ too.

We therefore may restrict to the case where $A$ is smooth and uniformly positive definite. 

\paragraph{The proof.}
We start as in the previous section, by writing
$$f=(f\det A)^{\frac1d}.$$
We apply Theorem 2.2 of \cite{YYL}~: given a matrix $S\in{\bf SDP}_d$ such that
\begin{equation}\label{detSA}
\det S=\Xint-_{\R^d/\Gamma} f(x)\,dx,
\end{equation}
there exists a $\Gamma$-periodic function $\phi_S\in{\mathcal C}^\infty$ such that $\det(S+\nabla^2\phi_S)=f$, and $S+\nabla^2\phi_S(x)\in{\bf SPD}_d$. In other words, the function $\psi_S(x)=\frac12\,x^TSx+\phi_S(x)$ is a convex solution of the Monge--Amp\`ere equation $\det\nabla^2\psi_S=f$.

Proceeding as in the bounded case, we have the inequality 
$$f\le\frac1d\,\Tr(A\nabla^2\psi_S)=\frac1d\,{\rm div}(A\nabla\psi_S)=\frac1d\,{\rm div}(A(Sx+\nabla\phi_S)).$$ 
Integrating over a fundamental domain, we obtain
\begin{equation}\label{RfA}
\Xint-_{\R^d/\Gamma} f(x)\,dx\le\frac1d\,\Xint-_{\R^d/\Gamma} (\Tr(AS)+{\rm div}(A\nabla\phi_S))\,dx=\frac1d\,\Tr(A_+S).
\end{equation}

There remains to minimize $\Tr(A_+S)$ under the constraint (\ref{detSA}). The minimum is achieved with $S=\lambda\,\widehat{A_+}$\,, where $\lambda$ is determined by
$$\lambda^d(\det A_+)^{d-1}=\Xint-_{\R^d/\Gamma} f(x)\,dx.$$
With this choice, (\ref{RfA}) becomes
$$\Xint-_{\R^d/\Gamma} f(x)\,dx\le\lambda\det A_+,$$
which is nothing but (\ref{Jdetper}).

\bigskip

The proof of Proposition \ref{p:caseeq} (the case of equality) is the same as that of Proposition \ref{p:eqcase}.

\bigskip

\paragraph{Proof of Corollary \protect\ref{c:CS}.}
Let $B(x):=A(x)-\bar A$, which satisfies ${\rm Div}\,B\equiv0$ and is compactly supported. Integrating by parts twice and using the assumption, we have
$$\int_\Omega b_{ij}dx = \int_\Omega b_{ij}\partial_ix_idx=-\int_\Omega x_i\partial_ib_{ij}dx  =\sum_{k\ne i}\int_\Omega x_i\partial_kb_{kj}dx=-\sum_{k\ne i}\int_\Omega b_{kj}\partial_kx_idx=0,
$$
whence the equality 
$$\Xint-_\Omega A(x)\,dx=\bar A.$$
Let $K$ be a cube containing $\Omega$. We extend $A$ to $K$ by setting $\hat A(x)=\bar A$ whenever $x\in K\setminus\Omega$. Next we extend $\hat A$ by periodicity to $\R^d$, $K$ being a fundamental domain. This $\hat A$ is a periodic DPT and has mean $\bar A$. Applying (\ref{Jdetper}) to $\hat A$, we have 
$$\int_\Omega(\det A)^{\frac1{d-1}}\,dx=\int_K(\det A)^{\frac1{d-1}}\,dx-(|K|-|\Omega|)(\det\bar A)^{\frac1{d-1}}\le(|K|-(|K|-|\Omega|))(\det\bar A)^{\frac1{d-1}},$$
from which the obtain the desired inequality.

\subsection{Proof of Theorem \protect\ref{th:nondiv}}

Let us establish first an {\em a priori} estimate when the tensor $A$ is smooth. We introduce as above the solution $\psi_S(x)=\frac12x^TSx+\phi_S(x)$ of the Monge--Amp\`ere equation $\det {\rm D}^2\psi_S=f:=(\det A)^{\frac1{d-1}}$, where $S\in{\bf SPD}_d$ is constrained by (\ref{detSA}) and $\phi_S$ is periodic. We still have
$$f=(f\det A)^{\frac1d}=(\det(A{\rm D}^2\psi_S))^{\frac1d}\le\frac1d\,\Tr(A{\rm D}^2\psi_S),$$
which writes now as
$$f\le\frac1d\,(\Tr(AS)+{\rm div}(A\nabla\phi_S)-({\rm Div}\,A)\nabla\phi_S).$$
Integrating over a fundamental domain, we obtain
$$\Xint-_{\R^d/\Gamma}f(x)\,dx\le\frac1d\,\Tr\left(S\Xint-_{\R^d/\Gamma}A(x)\,dx\right)+\frac1d\,\|{\rm Div}\,A\|_{\mathcal M}\sup_x|\nabla\phi_S(x)|.$$
To estimate the supremum of $\nabla\phi_S$, we involve the convexity of $\psi_S$. For every pair of points $x,x'$, we have
$$\psi_S(x')\ge\psi_S(x)+\nabla\psi_S(x)\cdot(x'-x),$$
that is
$$\frac12(x'-x)^TS(x'-x)+\phi_S(x')\ge\phi_S(x)+\nabla\phi_S(x)\cdot(x'-x),$$
When $x'-x=:\gamma$ is an element of the lattice $\Gamma$, the periodicity of $\phi_S$ yields
$$\frac12\gamma^TS\gamma\ge\nabla\phi_S(x)\cdot\gamma.$$
Replacing $\gamma$ by $-\gamma$, we actually have
\begin{equation}\label{phiSS}
|\nabla\phi_S(x)\cdot\gamma|\le\frac12\gamma^TS\gamma,\qquad\forall\gamma\in\Gamma\quad\hbox{and}\quad x\in\R^d.
\end{equation}
We now select a basis $(\gamma_1,\ldots,\gamma_d)$ of $\Gamma$, and form the matrix $M$ whose rows are the vectors $\gamma_j$. Writing $\nabla\phi_S=M^{-1}M\nabla\phi_S$ and using (\ref{phiSS}), we obtain the estimate
$$|\nabla\phi_S(x)|\le\|M^{-1}\|_{\ell^\infty\rightarrow\ell^2}\max_j\frac12\gamma_j^TS\gamma_j\le c_\Gamma\|S\|.$$
Because $S$ is non-negative, it satisfies $\|S\|\le\Tr\,S$. We therefore have
\begin{equation}\label{fDivS}
\Xint-_{\R^d/\Gamma}f(x)\,dx\le\frac1d\,\Tr S\left(\Xint-_{\R^d/\Gamma}A(x)\,dx+c_\Gamma\,\|{\rm Div}\,A\|_{\mathcal M}I_d\right).
\end{equation}
We choose as before 
$$S=\lambda \left(\Xint-_{\R^d/\Gamma}A(x)\,dx+c_\Gamma\,\|{\rm Div}\,A\|_{\mathcal M}I_d\right)$$
where $\lambda$ is determined by the condition (\ref{detSA}), and obtain the estimate
\begin{equation}\label{estnondiv}
\Xint-_{\R^d/\Gamma}(\det A(x))^{\frac1{d-1}}dx\le\left(\det \left(\Xint-_{\R^d/\Gamma}A(x)\,dx+c_\Gamma\,\|{\rm Div}\,A\|_{\mathcal M}I_d\right)\right)^{\frac d{d-1}}.
\end{equation}

\bigskip

We now turn to the general case.  Proceeding as in the proof of Theorem \ref{th:per}, we find a sequence of smooth positive definite symmetric tensors $A^\epsilon(x)=\epsilon I_d+\phi_\epsilon*A$, such that $A_\epsilon\rightarrow A$ in $L^1(\R^d/\Gamma)$. In addition, ${\rm Div}\,A^\epsilon=\phi_\epsilon*({\rm Div}\,A)$ converges vaguely towards ${\rm Div}\,A$. At last, up to an extraction, we may assume that $A_\epsilon(x)\rightarrow A(x)$ almost everywhere. We apply (\ref{estnondiv}) to each tensor $A_\epsilon$. Because of $\|{\rm Div}\,A_\epsilon\|_{\mathcal M}\le\|{\rm Div}\,A\|_{\mathcal M}$, we have
$$\Xint-_{\R^d/\Gamma}(\det A_\epsilon(x))^{\frac1{d-1}}dx\le\left(\det \left(\Xint-_{\R^d/\Gamma}A_\epsilon(x)\, dx+c_\Gamma\,\|{\rm Div}\,A\|_{\mathcal M}I_d\right)\right)^{\frac d{d-1}}.$$
We pass now to the limit as $\epsilon\rightarrow0+$. Because of 
$$\Xint-_{\R^d/\Gamma}A_{\epsilon}(x)\,dx\longrightarrow \Xint-_{\R^d/\Gamma}A(x)\,dx,$$ 
and using Fatou's Lemma, we recover (\ref{estnondiv}) for the tensor $A$. In particular, $(\det A(x))^{\frac1{d-1}}$ is integrable over the torus.

\section{Gas dynamics with finite mass and energy}\label{s:fluid}

\subsection{Proof of Theorem \protect\ref{th:absfl}}

Let us apply Theorem \ref{th:conv} in  the bounded convex domain $\Omega=(0,T)\times B_R$ for some $R>0$. We have
\begin{equation}\label{prelim}
\int_0^Tdt\int_{B_R}(\det A)^{\frac1n}dy\le \frac1{(n+1)|S^n|^{\frac1n}}\,\|A\vec n\|_{1,\partial\Omega}^{1+\frac1n}.
\end{equation}
The boundary consists in three parts: an initial ball $\{0\}\times B_R$, a final ball $\{T\}\times B_R$, and a lateral boundary $(0,T)\times S_R$. The latter contributes to 
$$g(R):=\int_0^Tdt\int_{S_R}\left|A\frac y{|y|}\right|\,dy.$$
Because $A$ is integrable, we have $g\in L^1(0,+\infty)$ and therefore there is a subsequence $R_m\rightarrow+\infty$ such that $g(R_m)\rightarrow0$. Passing to the limit in (\ref{prelim}), we obtain
\begin{equation}\label{nonhom}
\int_0^Tdt\int_{\R^n}(\det A)^{\frac1n}dy\le \frac1{(n+1)|S^n|^{\frac1n}}\,\left(\|(\rho, m)(0)\|_{L^1(\R^n)}+\|(\rho, m)(T)\|_{L^1(\R^n)}\right)^{1+\frac1n}.
\end{equation}
The latter estimate has the drawback that it is not homogeneous from a physical point of view. The density $\rho$ and the momentum $m$ have different dimensions and the norm
$$\|(\rho,m)\|_{L^1(\R^n)}=\int_{\R^n}\sqrt{\rho^2+|m|^2}\,dy$$
is not physically meaningful.

To recover the homogeneity, we introduce a scaling 
$$t'=\lambda t,\quad y'=y,\quad\rho'=\lambda^2\rho,\quad m'=\lambda m,\quad p'=p.$$
The corresponding $A'$ is a DPT over the slab $(0,T')\times\R^n$ where $T'=\lambda T$. Applying (\ref{nonhom}) to $A'$, we infer
$$\lambda^{1+\frac2n}\int_0^Tdt\int_{\R^n}(\det A)^{\frac1n}dy\le \frac1{(n+1)|S^n|^{\frac1n}}\,\left(\|(\lambda^2\rho,\lambda m)(0)\|_{L^1(\R^n)}+\|(\lambda^2\rho,\lambda m)(T)\|_{L^1(\R^n)}\right)^{1+\frac1n}.$$
Simplifying by $\lambda$ and then defining $\lambda=:\mu^{n+1}$, this becomes
\begin{eqnarray*}
\int_0^Tdt\int_{\R^n}(\det A)^{\frac1n}dy & \le & \frac1{(n+1)|S^n|^{\frac1n}}\,\left(\|(\mu^n\rho,\mu^{-1} m)(0)\|_{L^1(\R^n)}+\|(\mu^n\rho,\mu^{-1} m)(T)\|_{L^1(\R^n)}\right)^{1+\frac1n} \\
& \le & \frac1{(n+1)|S^n|^{\frac1n}}\,\left(2\mu^nM_0+\mu^{-1}(\|m(0)\|_{L^1(\R^n)}+\|m(T)\|_{L^1(\R^n))}\right)^{1+\frac1n}.
\end{eqnarray*}
We are free to choose the parameter $\mu>0$, and we make the choice
$$\lambda=\mu^{n+1}=\frac{\|m(0)\|_{L^1(\R^n)}+\|m(T)\|_{L^1(\R^n))}}{M_0}\,.$$
This yields the estimate in Theorem \ref{th:absfl}.

\subsection{The Euler and kinetic equations}

For the Euler equation, we only have to remark that $A$ is positive semi-definite and $\det A=\rho p^n$.

\bigskip

Likewise, for a kinetic equation, we only have to calculate the determinant of
$$A(t,y)=\begin{pmatrix} \int_{\R^n}f(t,y,v)\,dv & \int_{\R^n}f(t,y,v)v^T\,dv \\ \int_{\R^n}f(t,y,v)v\,dv & \int_{\R^n}f(t,y,v)v\otimes v\,dv \end{pmatrix}.$$
The formula
$$\begin{vmatrix} \int_{\R^n}f(v)\,dv & \int_{\R^n}f(v)v^T\,dv \\ \int_{\R^n}f(v)v\,dv & \int_{\R^n}f(v)v\otimes v\,dv \end{vmatrix}=\frac1{d!}\int\!\cdots\!\int_{(\R^n)^{n+1}}f(v^0)\cdots f(v^n)(\Delta(v^0,\ldots,v^n))^2dv^0\cdots dv^n$$
is a particular case of Andr\'eiev Identity
\begin{equation}\label{phiphiN}
\det\left(\!\!\left(\int\phi_i\phi_jd\mu(v)\right)\!\!\right)_{1\le i,j\le N}=\frac1{N!}\int^{\otimes N}\left(\det(\!(\phi_i(v_j)\,)\!)_{1\le i,j\le N}\right)^2d\mu(v_1)\cdots d\mu(v_N).
\end{equation}
To prove (\ref{phiphiN}), we develop the left-hand side as
$$\sum_{\sigma\in{\mathfrak S}_N}\epsilon(\sigma)\prod_i\int\phi_i(v)\phi_{\sigma(i)}(v)\,dv$$
and write 
$$\prod_i\int\phi_i(v)\phi_{\sigma(i)}(v)\,dv=\frac1{N!}\int^{\otimes N}\sum_{\rho\in{\mathfrak S}_N}\prod_i\phi_i(v^{\rho(i)})\phi_{\sigma(i)}(v^{\rho(i)})\,d\mu(v_1)\cdots d\mu(v_N).$$
There remains to verify
$$\sum_{\sigma\in{\mathfrak S}_N}\epsilon(\sigma)\sum_{\rho\in{\mathfrak S}_N}\prod_i\phi_i(v^{\rho(i)})\phi_{\sigma(i)}(v^{\rho(i)})=\left(\sum_{\lambda\in{\mathfrak S}_N}\epsilon(\lambda)\prod_i\phi_i(v^{\lambda(i)})\right)^2,$$
which is immediate.

\end{document}